\pdfoutput=1
\documentclass[a4paper,10 pt]{amsart}

\usepackage{color,soul}
\definecolor{red}{rgb}{1,0,0}
\setulcolor{red}

\usepackage{amssymb}
\usepackage{amsmath}
\usepackage{amsthm}
\usepackage{esint}
\usepackage{enumerate}

\usepackage[colorlinks=true,linkcolor=blue,citecolor=blue,urlcolor=blue]{hyper ref}
\usepackage{bookmark}

\swapnumbers
\newtheorem{lemma}{Lemma}[section]
\newtheorem{prop}[lemma]{Proposition}
\newtheorem{thm}[lemma]{Theorem}
\newtheorem{cor}[lemma]{Corollary}
\newtheorem{ann}[lemma]{Assumption}

\newtheorem*{Kap}{Theorem 1}

\theoremstyle{definition}
\newtheorem{defn}[lemma]{Definition}
\newtheorem{example}[lemma]{Example}
\newtheorem{rem}[lemma]{Remark}

\newtheorem{Satz}{Satz}[section]

\theoremstyle{definition}
\newtheorem{Definition}[Satz]{Definition}

\theoremstyle{remark}

\renewcommand{\(}{\left(}
\renewcommand{\)}{\right)}
\renewcommand{\~}{\tilde}
\renewcommand{\-}{\bar}

\newcommand{\cn}{\colon}

\newcommand{\N}{\mathbb{N}}
\newcommand{\R}{\mathbb{R}}

\renewcommand{\S}{\mathbb{S}}

\renewcommand{\a}{\alpha}
\renewcommand{\b}{\beta}
\newcommand{\g}{\gamma}
\renewcommand{\d}{\delta}
\newcommand{\e}{\epsilon}
\renewcommand{\k}{\kappa}

\newcommand{\vt}{\vartheta}
\newcommand{\s}{\sigma}
\newcommand{\p}{\varphi}
\newcommand{\G}{\Gamma}

\newcommand{\del}{\partial}

\DeclareMathOperator{\graph}{graph}

\DeclareMathOperator{\const}{const}
\DeclareMathOperator{\dist}{dist}
\DeclareMathOperator{\inte}{int}
\DeclareMathOperator{\Volu}{Vol}

\newcommand{\Def}{\begin{defn}}
\newcommand{\eDef}{\end{defn}}
\newcommand{\Thm}{\begin{thm}}
\newcommand{\eThm}{\end{thm}}
\newcommand{\Prop}{\begin{prop}}
\newcommand{\eProp}{\end{prop}}
\newcommand{\Rem}{\begin{rem}}
\newcommand{\eRem}{\end{rem}}
\newcommand{\Lem}{\begin{lemma}}
\newcommand{\eLem}{\end{lemma}}
\newcommand{\Ass}{\begin{ann}}
\newcommand{\eAss}{\end{ann}}
\newcommand{\Al}{\begin{align*}}
\newcommand{\eAl}{\end{align*}}
\newcommand{\eq}{\begin{equation}}
\newcommand{\eeq}{\end{equation}}
\newcommand{\ex}{\begin{example}}
\newcommand{\eex}{\end{example}}
\newcommand{\pf}{\begin{proof}}
\newcommand{\epf}{\end{proof}}
\newcommand{\Cor}{\begin{cor}}
\newcommand{\eCor}{\end{cor}}

\newcommand{\Ra}{\Rightarrow}
\newcommand{\ra}{\rightarrow}

\newcommand{\mc}{\mathcal}

\newcommand{\mrm}{\mathrm}

\newcommand{\hp}{\hphantom}

\newcommand{\bfi}{\begin{fig}}
\newcommand{\efi}{\end{fig}}
\newcommand{\btab}{\begin{tab}}
\newcommand{\etab}{\end{tab}}
\newcommand{\barr}{\begin{array}}
\newcommand{\earr}{\end{array}}
\newcommand{\beqq}{\begin{equation}}
\newcommand{\eeqq}{\end{equation}}
\newcommand{\beao}{\begin{align*}}
\newcommand{\eeao}{\end{align*}\noindent}
\newcommand{\beam}{\begin{eqnarray}}
\newcommand{\eeam}{\end{eqnarray}\noindent}
\newcommand{\bdis}{\begin{displaymath}}
\newcommand{\edis}{\end{displaymath}\noindent}

\newcommand{\bbn}{\mathbb{N}}

\newcommand{\bbr}{\mathbb{R}}

\newcommand{\bbs}{\mathbb{S}}

\newcommand{\cals}{{\mathcal S}}

\newcommand{\calh}{{\mathcal H}}
\newcommand{\calp}{\mathcal{P}}

\newcommand{\dt}{\frac{d}{dt}}

\newcommand{\ov}{\overline}

\DeclareMathOperator{\interior}{int}

\begin{document}

\numberwithin{equation}{section}

\title[Inverse curvature flows in the sphere]{Rigidity results, inverse curvature flows and Alexandrov-Fenchel type inequalities in the sphere}
\author{Matthias Makowski and Julian Scheuer}

\thanks{The second author has been supported by the DFG}
\date{\today}

\address{Dr. Matthias Makowski, Universit\"at Konstanz, 
Fachbereich Mathematik und Statistik,
78457 Konstanz, Germany}
\email{matthias.makowski@uni-konstanz.de}

\address{Dr. Julian Scheuer\ \\ Ruprecht-Karls-Universit\"at, Institut f\"ur Angewandte Mathematik, Im Neuenheimer Feld 294, 69120 Heidelberg, Germany}
\email{scheuer@math.uni-heidelberg.de}

\begin{abstract}
We prove a rigidity result in the sphere which allows us to generalize a result about smooth convex hypersurfaces in the sphere by Do Carmo-Warner to convex $C^{2}$-hypersurfaces. We apply these results to prove $C^{1,\beta}$-convergence of inverse $F$-curvature flows in the sphere to an equator in $\S^{n+1}$ for embedded, closed, strictly convex initial hypersurfaces. The result holds for large classes of curvature functions including the mean curvature and arbitrary powers of the Gauss curvature. We use this result to prove Alexandrov-Fenchel type inequalities in the sphere.
\end{abstract}

\maketitle
\makeatletter
\def\l@subsection{\@tocline{2}{0pt}{2.5pc}{5pc}{}} 
\makeatother

\tableofcontents

\section{Introduction}
\label{Introduction}

This work deals with geometric problems on the $(n+1)$-dimensional unit sphere $\S^{n+1}\subset\{x\in \R^{n+2}\cn |x|=1\}.$ We assume $n\geq 2,$ unless stated otherwise.
We are interested in the connection between (analytically) convex hypersurfaces and (geodesically) convex bodies. The notion of convexity of sets is significantly more subtle than in Euclidean space, due to the 
existence of focal points in the sphere. \\
A very well known result in this direction by do Carmo and Warner, \cite[Theorem 1.1]{DoCarmoWarner1970}, is the following.
\begin{Kap}[Do Carmo, Warner]
\label{doCarmoWarner}
Let $x\cn M^{n}\ra\S^{n+1}$ be an isometric immersion of a compact, connected, orientable $n$-dimensional $C^{\infty}$-Riemannian manifold into the $(n+1)$-sphere of sectional curvature equal to one, and assume that all sectional curvatures of $M^{n}$ are greater than or equal to one. Then $x$ is an embedding, $M^{n}$ is diffeomorphic with $\S^{n}$ and $x(M^{n})$ is either totally geodesic or contained in an open hemisphere. In the latter case $x(M^{n})$ is the boundary of a convex body in $\S^{n+1}.$
\end{Kap}

Also compare \cite{Ando}, which deals with strictly convex hypersurfaces.

 We will show that some parts of this result can be generalized to non-smooth, geodesically convex bodies in the sphere. In particular we will prove the following result, for the exact definitions of weakly convex bodies see Section \ref{RigiditySection}.
 
 \Thm
\label{RigidityResult}
	Let $ n\geq 1$ and $\hat{M}\subset \bbs^{n+1}$ be a weakly convex body in a hemisphere. Let $x_0\in \bbs^{n+1}$ be such that $\hat{M}$ is contained in the closed hemisphere $\calh(x_0)$ with equator $\cals(x_{0})$. Suppose that $\hat{M}$ satisfies an interior sphere condition at all points $p\in \hat{M} \cap \cals(x_0)$. Then either $\hat{M}$ is equal to $\calh(x_0)$ or $\hat{M}$ is contained in an open hemisphere. 
\eThm

With the help of this result, we prove that the strong regularity assumption in \cite{DoCarmoWarner1970} is not necessary.

\Cor
\label{C2Rigidity}
	Let $M\subset \bbs^{n+1}$ be an embedded, closed, connected and convex $C^{2}$-hypersurface. Then $M$ is either an equator or $M$ is contained in an open hemisphere and bounds a convex body.
\eCor

We apply those rigidity results to treat an inverse curvature flow in the sphere $\S^{n+1}$ of the form

\begin{equation}
\label{floweq}
\begin{split}
	\dot{x} &= -\Phi(F)\, \nu,\\ 
	x(0) &= x_0,
\end{split}
\end{equation}
where $x_0: \bbs^n \rightarrow \bbs^{n+1}$ is the embedding of an initial hypersurface $M_0 := x_0(\bbs^n)$ of class $C^{4, \alpha}$ for some $0<\alpha <1$, which is furthermore required to be strictly convex. $\nu$ is the corresponding outer normal, $\Phi\in C^\infty(\bbr_+, \bbr)$, $\Phi(x) = -x^{-p}$, $p>0,$ $F$ is a curvature function evaluated at the principal curvatures of the flow hypersurfaces $M_t$ and $x(t)$ denotes the embedding of $M_t$. 
We will show that under certain assumptions, cf. \ref{MainAssumption}, the flow exists up to a finite time and converges in $C^{1,\b}$ to the embedding of an equator.

Curvature flows and their application to geometric inequalities have been treated for over thirty years. Following the ground breaking work by G. Huisken, \cite{gh:mean}, who considered the mean curvature flow, also inverse, or expanding flows have been considered. Here the works on the inverse curvature flow by C. Gerhardt, \cite{cg90}, as well as J. Urbas, \cite{Urbas}, have to be mentioned, where also non-convex hypersurfaces were considered.
Similar results have been shown in other ambient spaces and for general $p$-homogeneous curvature functions, e.g. \cite{cg:ImcfHyp}, \cite{cg:Impf}, or \cite{js:ImpfHyp}.

We consider a large class of curvature functions. We allow other homogeneities than $1,$ in particular our result holds for arbitrary powers of the Gaussian curvature without further pinching assumptions. The detailed assumptions on the curvature function are listed below, whereafter we state the convergence result.

\Ass
\label{MainAssumption}
Suppose $F\in C^{2,\a}(\G),$ $0<\a<1,$ is a symmetric function, where $\Gamma$ is the positive cone $\Gamma_+= \{\kappa = (\kappa_i)\in \bbr^n\cn \kappa_i > 0 \,\, \forall \, i \in \{1, \ldots, n\}\}$. We need the following assumptions for the curvature function $F$:
\begin{itemize}
\item $F$ is positively homogeneous of degree $1$, i.e. $\forall\, \kappa \in \Gamma_+$, $\forall\, \lambda \in \bbr_+$: $F(\lambda \kappa) = \lambda F(\kappa)$.
\item $F$ is strictly increasing in each argument: $\forall\, i \in \{1, \ldots, n\}$, $\forall\, \kappa \in \Gamma_+$ there holds $F_i(\kappa)$ = $\frac{\partial F}{\partial \kappa_i}(\kappa) > 0$.
\item $F$ is positive, $F_{|\Gamma_+} > 0$, and $F$ is normalized, $F(1, \ldots, 1) = n$.
\item Either:
\begin{enumerate}[(i)]
\item $F$ is concave and inverse concave, i.e. $F_{-1}(\kappa_i) := \frac{1}{F(\kappa_i^{-1})}$ is concave. 
\item $F$ is concave and $F$ approaches zero on the boundary of $\Gamma_+$.
\end{enumerate}
\item If $p\neq1$, we assume (ii) is valid.
\end{itemize} 
\eAss


The most important examples of curvature functions $F$ being concave and inverse concave are $\left(\frac{H_k}{H_l}\right)^{\frac{1}{k-l}}$, $n\geq k> l \geq 0$, or the power means $\left(\sum_{i=1}^n \kappa_i^r \right)^{\frac{1}{r}}$ for $|r|\leq 1$ . For a proof of the inverse concavity of these functions see the proofs of \cite[Theorem 2.6, Theorem 2.7]{AndrewsPinching}. 
Our exact result concerning the curvature flows is:

\begin{thm} \label{mainthm1}

Let $0<\alpha<1$. Let $\S^{n}\hookrightarrow M_{0}\subset \S^{n+1}$ be an embedded,  strictly convex hypersurface of class $C^{4,\a}.$ Let $F$ be a curvature function satisfying \ref{MainAssumption}. Then there exists a finite time $0<T^*<\infty$ and a unique curvature flow $$x\in H^{2+\a,\frac{2+\a}{2}}( [0,T^{*})\times \S^{n},\S^{n+1}),$$ 
which satisfies the flow equation
\begin{align} \begin{split}
\dot{x}&=F^{-p}\nu\\
						x(0)&=M_{0},\end{split}\end{align}
	where $0<p<\infty,$ $\nu(t,\xi)$ is the outward normal to $M_{t}=x(t,M)$ at $x(t,\xi)$ and there exists $0<t_0<T^*$ such that the leaves $M_{t}$, $t_0\leq t<T^*$, are graphs over some suitable equator $\cals({x_{0}}),$ $x_{0}\in\S^{n+1},$
	\eq M_{t}=\graph u(t,\cdot),\eeq
	where $u$ is the radial distance to $x_{0}.$
For $t\ra T^{*},$ the functions $u(t,\cdot)$ converge to $\frac{\pi}{2}$ in $C^{1,\beta}(\S^{n})$ for arbitrary $0<\beta<1$ and we have for $1\leq q<\infty,$ that
\eq \int_{M_{t}}H^{q}\ra 0,\ t\ra T^{*}.\eeq 
\end{thm}

In this theorem, $H^{2+\a,\frac{2+\a}{2}}( [0,T^{*})\times \S^{n},\S^{n+1})$ denotes the parabolic Hoelder space as in \cite[Definition 2.5.2]{cg:cp}.

Recently, Gerhardt also considered inverse curvature flows of strictly convex hypersurfaces in $\bbs^{n+1}$ by curvature functions satisfying the assumptions of \ref{MainAssumption}(i) in case $p=1,$ see \cite{Gerhardt:/2015}. He obtains smooth convergence of the flow to an equator. However, his methods substantially differ from ours.

Theorem \ref{mainthm1} allows us to prove Alexandrov-Fenchel type inequalities in the sphere, namely:

\Thm
Let $M\subset \bbs^{n+1}$ be an embedded, closed, connected and convex $C^2$-hypersurface of the sphere. Then we have the inequality
	\begin{equation}
		(\tilde{V}_1(M))^2 \geq (\tilde{V}_0(M))^{2\left(\frac{n-1}{n}\right)} - (\tilde{V}_0(M))^2,
	\end{equation}
	and equality holds if and only if $M$ is a geodesic sphere.
	
	Furthermore, if $n\geq 3,$ we have the inequality
	\begin{equation}
		\tilde{V}_2(M) \geq \(\tilde{V}_0(M)\)^{\frac{n-2}{n}} - \tilde{V}_0(M),
	\end{equation}
	and equality holds if and only if $M$ is a geodesic sphere.
	
	Let $k\in \bbn_+$ with $2k+1\leq n$ and let $\hat{M}$ be the convex body enclosed by $M$. Then we have the inequality
	\begin{equation}
	\label{GeomIneqEqIII}
		W_{2k+1}(\hat{M}) \geq \frac{\omega_n}{n+1}Ê\sum_{i=0}^k(-1)^i\frac{n-2k}{n-2k+2i} {k\choose i} \left(\frac{n+1}{\omega_n}W_1(\hat{M})\right)^{\frac{n-2k+2i}{n}}.
	\end{equation}	
	and equality holds if and only if $M$ is a geodesic sphere.
\eThm
Here $\tilde{V}_k(M)$ denotes, up to a constant, the $k$-th mean curvature integral and is defined by
\begin{equation}
	\tilde{V}_k(M) := \omega_n^{-1}\int_{M} \tilde{H}_k \,d\mu,
\end{equation}
where $\tilde{H}_k := \frac{H_k}{{n\choose k}}$ are the normalized elementary symmetric polynomials and $\omega_n := |\bbs^n|$. $W_k(\hat{M})$ denotes the $k$-th quermassintegral of $\hat{M}$, see Section 7 for a definition. For a more detailed account of the mean curvature integrals and their relation to the quermassintegrals in spaces of constant curvature, see for example \cite{SolanesInt}. Especially inequality \eqref{GeomIneqEqIII} shows that the geometric inequalities for the quermassintegrals resemble the corresponding inequalities in hyperbolic space, see \cite[Theorem 1.3]{GeWangWuII}.

Curvature flows have shown to be a useful method to obtain geometric inequalities. Probably the most known result in this direction is the proof of the Riemannian Penrose Inequality by Huisken and Ilmanen in \cite{HuisIlm} using an inverse mean curvature flow in asymptotically flat $3$-manifolds. 

But also Alexandrov-Fenchel type inequalities have been proved using curvature flows:
In Euclidean space, McCoy showed in \cite{McCoyMixedAreaGen}, that the Alexandrov-Fenchel inequalities for strictly convex hypersurfaces can be deduced from 
a mixed-volume preserving curvature flow. In 2009, Guan and Li (see \cite{GuanLi}) used inverse $F$-curvature flows in Euclidean space to show these inequalities for $k$-convex, starshaped domains. Recently, the first author transferred the results about mixed-volume preserving curvature flows in Euclidean space from \cite{McCoyMixedAreaGen} to the hyperbolic space in \cite{MakowskiHypVol}. Wang and Xia used these results in \cite{WangXia} to obtain the Alexandrov-Fenchel inequalities for horospherically convex hypersurfaces in hyperbolic space. Some of these inequalities have also been shown in hyperbolic space by using inverse $F$-curvature flows, see for example \cite{GeWangWuI}, \cite{GeWangWuII}, \cite{LimaGirao}.

\section{Setting and general facts}
We now state some general facts about hypersurfaces, especially those that can be written as graphs. We basically follow the description of \cite{cg:ImcfHyp} and \cite{js:ImpfHyp}, but restrict to Riemannian manifolds. For a detailed discussion we refer to \cite{cg:cp}. \ \\
Let $N=N^{n+1}$ be Riemannian and $M=M^{n}\hookrightarrow N$ be a hypersurface. The geometric quantities of $N$ will be denoted by $(\-{g}_{\a\b}),$ $(\-{R}_{\a\b\g\d})$ etc., where greek indices range from $0$ to $n$. Coordinate systems in $N$ will be denoted by $(x^{\a}).$
Quantities for $M$ will be denoted by $(g_{ij}),$ $(h_{ij})$ etc., where latin indices range from $1$ to $n$ and coordinate systems will generally be denoted by $(\xi^{i}),$ unless stated otherwise.\ \\
Covariant differentiation will usually be denoted by indices, e.g. $u_{ij}$ for a function $u\colon M\rightarrow \R$, or, if ambiguities are possible, by a semicolon, e.g.
$h_{ij;k}.$ Usual partial derivatives will be denoted by a comma, e.g. $u_{i,j}.$\ \\
Let $x\colon M\hookrightarrow N$ be an embedding and $(h_{ij})$ be the second fundamental form with respect to a normal $-\nu,$ i.e. we have the \textit{Gaussian formula}
\eq x^{\a}_{ij}=-h_{ij}\nu^{\a}, \eeq where $\nu$ is a differentiable normal, the \textit{Weingarten equation}
\eq \nu^{\a}_{i}=h^{k}_{i}x^{\a}_{k}, \eeq
the \textit{Codazzi equation} \eq h_{ij;k}-h_{ik;j}=\-{R}_{\a\b\g\d}\nu^{\a}x^{\b}_{i}x^{\g}_{j}x^{\d}_{k} \eeq
and the \textit{Gau\ss\ equation} \eq R_{ijkl}=(h_{ik}h_{jl}-h_{il}h_{jk})+\-{R}_{\a\b\g\d}x^{\a}_{i}x^{\b}_{j}x^{\g}_{k}x^{\d}_{l}. \eeq

Now assume that $N=(a,b)\times S_{0},$ where $S_{0}$ is compact Riemannian and that there is a Gaussian coordinate system $(x^{\a})$ such that
\eq d\-{s}^{2}=e^{2\psi}((dx^{0})^{2}+\sigma_{ij}(x^{0},x)dx^{i}dx^{j}), \eeq
where $\sigma_{ij}$ is a Riemannian metric, $x=(x^{i})$ are local coordinates for $\mc{S}_{0}$ and $\psi\colon N\rightarrow \R$ is a function.\ \\
Let $M=\graph u_{|\mc{S}_{0}}$ be a hypersurface
\eq M=\{(x^{0},x)\colon x^{0}=u(x), x\in \mc{S}_{0}\}, \eeq
then the induced metric has the form \eq g_{ij}=e^{2\psi}(u_{i}u_{j}+\sigma_{ij}) \eeq
with inverse \eq g^{ij}=e^{-2\psi}(\sigma^{ij}-v^{-2}u^{i}u^{j}), \eeq
where $(\sigma^{ij})=(\sigma_{ij})^{-1},$ $u^{i}=\sigma^{ij}u_{j}$ and \eq v^{2}=1+\sigma^{ij}u_{i}u_{j}\equiv 1+|Du|^{2}.\eeq
We use, especially in the Gaussian formula, the normal 
\eq \label{outernormal} (\nu^{\a})=v^{-1}e^{-\psi}(1,-u^{i}). \eeq
Looking at $\a=0$ in the Gaussian formula, we obtain 
\eq e^{-\psi}v^{-1}h_{ij}=-u_{ij}-\-{\Gamma}^{0}_{00}u_{i}u_{j}-\-{\Gamma}^{0}_{0i}u_{j}-\-{\Gamma}^{0}_{0j}u_{i}-\-{\Gamma}^{0}_{ij} \eeq
and \eq e^{-\psi}\-{h}_{ij}=-\-{\Gamma}^{0}_{ij}, \eeq
where covariant derivatives are taken with respect to $g_{ij}.$

In our special situation $N=\S^{n+1}$ let $x_{0}\in\S^{n+1},$ then  by introducing geodesic polar coordinates we derive a representation of the metric in the form
\eq d\-{s}^{2}=dr^{2}+\sin^{2}r\s_{ij}dx^{i}dx^{j},\eeq 
where $\s_{ij}$ is the canonical metric of $\S^{n}$ and $0<r<\pi.$ Then we obtain for a geodesic sphere given by a constant graph $u\equiv r$ with $0<r<\pi$ that $\bar{h}_{ij} = \frac{\bar{H}}{n} \bar{g}_{ij}$ and 
\begin{equation}
\label{HLevelHyp}
	\frac{\bar{H}}{n}(r) = \frac{\cos r}{\sin r}. 
\end{equation}

 Using \cite[Theorem 1.1]{DoCarmoWarner1970}, we conclude for an embedding of a smooth, strictly convex, closed hypersurface $M,$ that it is contained in an open hemisphere and thus it can be written as a graph over $\S^{n}$ in the previously described coordinate system, i.e.
\eq M=\graph u_{|\mc{S}_{0}}.\eeq

Now we want to give some elementary facts about the curvature functions. Firstly, we provide the definition of these functions and mention some identifications, which will be used in the sequel without explicitly stating them again.
\begin{Definition}
	Let $\Gamma \subset \bbr^n$ be an open, convex, symmetric cone, i.e.
	\begin{equation}
		(\kappa_i) \in \Gamma \Longrightarrow (\kappa_{\pi i}) \in \Gamma \quad \forall \, \pi \in \mathcal{P}_n,
	\end{equation}
	where $\mathcal{P}_n$ is the set of all permutations of order $n$. Let $f \in C^{m, \alpha}(\Gamma)$, $m \in \bbn$, $0\leq \a \leq 1$, be \textit{symmetric}, i.e.,
	\begin{equation}
		f(\kappa_i) = f(\kappa_{\pi i}) \quad \forall \, \pi \in \mathcal{P}_n.
	\end{equation}
	Then $f$ is said to be a \textit{curvature function} of class $C^{m,\alpha}$. For simplicity we will also refer to the pair $(f, \Gamma)$ as a curvature function.
\end{Definition}
Now denote by $\mathbf{S}$ the symmetric endomorphisms of $\bbr^n$ and by $\mathbf{S}_\Gamma$ the symmetric endomorphisms with eigenvalues belonging to $\Gamma$, an open subset of $\mathbf{S}$. If $(f, \Gamma)$ is a smooth curvature function, we can define a mapping
\begin{equation}
\begin{split}
	F: &\,\mathbf{S}_\Gamma \rightarrow \bbr,\\
	&A\mapsto f(\kappa_i),
\end{split}
\end{equation}
where the $\kappa_i$ denote the eigenvalues of $A$. For the relation between these different notions, especially the differentiability properties and the relation between their derivatives, see \cite[Chapter 2.1]{cg:cp}. Since the differentiability properties are the same for $f$ as for $F$ in our setting, see \cite[Theorem 2.1.20]{cg:cp}, we do not distinguish between these notions and always write $F$ for the curvature function. Hence at a point $x$ of a hypersurface we can consider a curvature function $F$ as a function defined on a cone $\Gamma \subset \bbr^n$, $F = F(\kappa_i)$ for $(\kappa_i) \in \Gamma$ (representing the principal curvatures at the point $x$ of the hypersurface), as a function depending on $(h_i^j)$, $F = F(h_i^j),$ or as a function depending on $(h_{ij})$ and $(g_{ij})$, $F = F(h_{ij}, g_{ij})$. However, we distinguish between the derivatives with respect to $\Gamma$ or $\mathbf{S}$. We briefly summarize our notation and important properties:

	For a smooth curvature function $F$ we denote by $F^{ij} = \frac{\partial F}{\partial h_{ij}}$, a contravariant tensor of order 2, and $F^j_i = \frac{\partial F}{\partial h_j^i}$, a mixed tensor, contravariant with respect to the index $j$ and covariant with respect to $i$. We also distinguish the partial derivative $F_{,i} = \frac{\partial F}{\partial \kappa_i}$ and the covariant derivative $F_{;i} = F^{kl}h_{kl;i}$. Furthermore $F^{ij}$ is diagonal if $h_{ij}$ is diagonal and in such a coordinate system there holds $F^{ii} = \frac{\partial F}{\partial \kappa_i}$. For a relation between the second derivatives see \cite[Lemma 2.1.14]{cg:cp}. Finally, if $F \in C^2(\Gamma)$ is concave (convex), then $F$ is also concave (convex) as a curvature function depending on $(h_{ij})$. 


\section{Rigidity results}
\label{RigiditySection}
In this section we want to prove the rigidity result, Theorem \ref{RigidityResult} and Corollary \ref{C2Rigidity}. First we need some definitions, which also apply to the case $n\geq 1.$

\Def

\begin{enumerate}[(i)]
\item For a point $x\in \bbs^{n+1}$ we will denote the closed hemisphere with center in $x$ by $\calh(x)$, 
\begin{equation}
	\calh(x) := \{p\in \bbs^{n+1}: \langle p, x\rangle \geq 0\},
	\end{equation}
	where $\langle\cdot,\cdot\rangle$ is the scalar product in $\R^{n+2},$
and the corresponding equator by $\cals(x) := \calh(x)\setminus\interior\calh(x)$.

\item For points $p, q \in \bbs^{n+1},$ $\hat{\gamma}_{p,q}:[0, L]\rightarrow \bbs^{n+1}$ will denote a $C^1$ parametrization by arc length of the geodesic segment $\hat\Gamma_{p, q}$ with $\hat\gamma_{p,q}(0) = p$ and $\hat\gamma_{p,q}(L) = q$. The geodesic segment is not unique if $\dist(p, q) = \pi$. 

\item For $x\in \bbs^{n+1}$ the stereographic projection mapping $x$ to $0$ will be denoted by $\calp_x:\bbs^{n+1}\setminus\{-x\}\rightarrow \bbr^{n+1}$.
\end{enumerate}
\eDef

\Def

\begin{enumerate}[(i)] 
\item Let $\hat{M} \subset \bbs^{n+1}$ be a set. We say that $\hat{M}$ is a weakly convex set in $\bbs^{n+1}$ [in a hemisphere], if [there exists $x \in \bbs^{n+1}$ such that $\hat{M} \subset \calh(x)$ and] for arbitrary $p, q \in \hat{M}$ there exists a minimizing geodesic $\hat\Gamma_{p,q}$ connecting $p$ and $q$, which is contained in $\hat{M}$. 
\item Let $\hat{M} \subset \bbs^{n+1}$ be a set. We say that $\hat{M}$ is a convex set in $\bbs^{n+1}$ [in a hemisphere], if [there exists $x \in \bbs^{n+1}$ such that $\hat{M} \subset \calh(x)$ and] for arbitrary $p, q \in \hat{M}$ all minimizing geodesics $\hat\Gamma_{p,q}$ connecting $p$ and $q$ [and contained in $\calh(x)$] are contained in $\hat{M}$.
\item We say that $\hat{M} \subset \bbs^{n+1}$ is a (weakly) convex body [in a hemisphere], if it is a compact, (weakly) convex set [in a hemisphere] with nonempty interior.
\item A set $MÊ\subset \bbs^{n+1}$ is a closed, geodesically convex hypersurface, if there exists a convex body $\hat{M} \subset\bbs^{n+1}$ in a hemisphere, such that $M = \partial \hat{M}$. The set $\hat{M}$ is called the convex body of $M$.
\item Let $x\in \bbs^{n+1}$. Let $\hat{M}\subset \bbs^{n+1}$ be a set with $\ov{\hat{M}} \subset \bbs^{n+1}\setminus\{-x\}$. We say that $\hat{M}$ satisfies an interior sphere condition at a point $p\in \ov{\hat{M}}$ with respect to $x$, if the set $\calp_x(\ov{\hat{M}})\subset \bbr^{n+1}$ satisfies an interior sphere condition at $\calp_x(p)$. 
\end{enumerate}
\eDef
\Rem
	The following observations have to be made:
	\begin{enumerate}[(i)]
		\item Let $M\subset \bbs^{n+1}$ be a closed, geodesically convex hypersurface. Then the convex body of $M$ is not unique, as can be seen by looking at $M= \cals(x)$, where $x\in \bbs^{n+1}$ is arbitrary.
		\item Note that there are different notions of convexity in the sphere: we do not demand that a geodesic connecting two points in the convex body has to be unique. It is well-known, see also Lemma \ref{OpenHemisphereLemma}, that a convex body $\hat{M}$ in the sphere, which does not contain a pair of antipodal points is contained in an open hemisphere. If on the other hand the convex body $\hat{M}$ in the sphere contains antipodal points, then it follows from the definition that $\hat{M} = \bbs^{n+1}$. 
	\end{enumerate}
\eRem

\Rem \label{StereographicProjection}
The following observations can be found in \cite[p. 278, 279]{cg:cp}. Let $x_0 \in \bbs^{n+1}$. Defining 
\eq \rho=2\tan\frac{r}{2},\eeq
where $r$ is the geodesic distance to $x_{0}$,
we obtain a representation of the spherical metric as
\eq d\-{s}^{2}=\frac{1}{\(1+\frac{1}{4}\rho^{2}\)^{2}}(d\rho^{2}+\rho^{2}\s_{ij}d\xi^{i}d\xi^{j})\equiv e^{2\psi}\hat{g}_{\a\b},\eeq
where $\hat{g}$ denotes the Euclidean metric in $\bbr^{n+1}$. A point $q \in \bbs^{n+1}$ is contained in $\calh(x_0)$ if and only if $r\leq \frac{\pi}{2}$, which is equivalent to $\rho\leq 2$. 

A $C^2$-hypersurface $M\subset \bbs^{n+1}\setminus\{-x_0\}$ can be seen as embedded in Euclidean space using the conformally flat parametrization of $\S^{n+1}$ via stereographic projection $\calp_{x_0}$. We will denote the hypersurface $\calp_{x_0}(M)$ by $\mc{M}$. The second fundamental form $h^i_j$ of $M$ and the corresponding Euclidean quantity $\hat{h}^{i}_{j}$ are related by
\eq \label{Ster.Proj.A} e^{\psi}h^{i}_{j}=\hat{h}^{i}_{j}+\psi_{\a}\hat{\nu}^{\a}\d^{i}_{j},\eeq
where $\hat\nu$ denotes the Euclidean normal vector field of $\mc{M}$.
Thus a simple calculation reveals that for a strictly convex and $C^{2}$-bounded $M$, the corresponding hypersurface $\mc{M}$ is strictly convex and bounded in $C^{2}.$
\eRem

The closure of a weakly convex set is again a weakly convex set. However, this statement is not true for convex sets (neither in spheres nor in hemispheres). We want to prove a sufficient condition for a weakly convex body in a hemisphere to be a convex body in a hemisphere.
\Thm
\label{WeaklyConvexTheorem}
	Let $n\geq 1$ and $\hat{M}\subset \bbs^{n+1}$ be a weakly convex body in a hemisphere $\calh(x_0)$ for some $x_0 \in \bbs^{n+1}$. Suppose that $\hat{M}$ satisfies an interior sphere condition with respect to $x_0$ at all points $p\in \hat{M} \cap \cals(x_0)$. Then $\hat{M}$ is a convex body in a hemisphere.
\eThm

Firstly, we need some lemmata.

\Lem
\label{CharacterizationLemma}
	Let $n\geq 1$ and $x_0\in \bbs^{n+1}$. Let $p\in \cals(x_0)$. Let $\gamma:[0, \pi]\rightarrow \bbs^{n+1}$ be a $C^1$-geodesic, parametrized by arc length, with $\gamma(0) = p$. Let $\tilde{p}$ denote the outward normal vector of $\calh(x_0)$ at $p$. Then 
	\begin{equation}
		\langle\dot\gamma(0), \tilde{p}\rangle < (>, =)\,,
	\end{equation}
	if and only if the geodesic satisfies $\gamma(t) \in \interior\calh(x_0)\, \(\complement \calh(x_0), \cals(x_0)\)$ for some (and hence every) $t \in (0, \pi)$.
\eLem
\begin{proof}
	First of all, we note the following fact: Since $\gamma$ is a $C^1$-geodesic, $\gamma$ is a segment of a great circle. Hence a third point lying on $\gamma$ determines uniquely the great circle $\Gamma$ such that $\gamma([0, \pi]) \subset \Gamma$. Thus the existence of $t\in (0, \pi)$, such that $\gamma(t) \in \interior\calh(x_0)\, (\complement \calh(x_0), \cals(x_0))$ implies $\gamma((0, \pi)) \subset  \interior\calh(x_0)\, (\complement \calh(x_0), \cals(x_0))$.

	Suppose firstly, that $\langle \dot\gamma(0), \tilde{p}\rangle < 0$. Then $\gamma(t) \in \interior\calh(x_0)$ for $t$ close to $0$. Hence by the observation made above, we obtain $\gamma((0, \pi)) \subset  \interior\calh(x_0)$. 

	If on the other hand there exists $t \in (0, \pi)$, such that $\gamma(t) \in \interior\calh(x_0)$, then the geodesic $\tilde{\gamma}:[0, t] \rightarrow \bbs^{n+1}$ with $\tilde{\gamma}(0) = \gamma(t)$ and $\dot{\tilde\gamma}(s) = -\dot\gamma(t-s)$ satisfies $\langle \dot{\tilde\gamma}(t), \tilde{p}\rangle > 0$, hence we obtain $\langle \dot\gamma(0), \tilde{p}\rangle < 0$.
\end{proof}


\Lem
\label{StereographicInball}
	Let $n\geq 1$ and $\hat{M}\subset \bbs^{n+1}$ be a weakly convex body in the hemisphere $\calh(x_0)$ for some fixed $x_0\in \bbs^{n+1}$.
	Let $p\in \hat{M} \cap \cals(x_0)$ and suppose that $\hat{M}$ satisfies an interior sphere condition at $p$ with respect to $x_0$. Let $\gamma \in C^1([0, t_0), \bbs^{n+1})$, $0<t_0$, be a geodesic with $\gamma(0) = p$ and $\gamma((0, t_0)) \in \interior \calh(x_0)$. Then there exists $0<\delta\leq t_0$, such that $\gamma((0, \delta)) \subset \interior \hat{M}$.
\eLem
\begin{proof}
	Let $\tilde{p}$ denote the outward normal vector of $\calh(x_0)$ at $p$.	
	From Lemma \ref{CharacterizationLemma} we obtain 
	\begin{equation}
	\label{NormalInequality}
		\langle \dot\gamma(0), \tilde{p}\rangle < 0.
	\end{equation}
	
	Let us look at the situation in the coordinates given by the stereographic projection $\calp_{x_0}.$	
	        	Let $\rho>0$ and $B_\rho(\bar{p})$ be an inball with respect to $\calp_{x_0}(p)$ with center $\bar{p}.$		Let $\tilde{\gamma} := \calp_{x_0}\circ\gamma$, then since the metric of the sphere is conformally equivalent to the Euclidean metric in stereographic coordinates, we obtain from \eqref{NormalInequality}
	\begin{equation}
		 \left\langle \dot{\tilde\gamma}(0), \nu(\calp_{x_0}(p))\right\rangle < 0,
	\end{equation}
	where $\nu(\calp_{x_0}(p))$ denotes the outward normal of $B_2(0)$ at $\calp_{x_0}(p).$
	
	Since the inball $B_{\rho}(\bar{p})$ is tangent to $\partial B_2(0)$, we obtain some small $\delta > 0$, such that for $t\in (0, \delta)$:
\begin{equation}
	\tilde\gamma(t) \in B_\rho(\bar{p}).
\end{equation}
\end{proof}


\begin{proof}[Proof of Theorem \ref{WeaklyConvexTheorem}]
	Let $p, q\in \hat{M}$ be two arbitrary points, then we have to show that an arbitrary minimizing geodesic $\hat\gamma_{p,q}$ connecting $p$ and $q$ and contained in $\calh(x_0)$ is contained in $\hat{M}$.

If no antipodal points exist in $\hat{M}$, we have nothing to prove. Hence let us assume that there exist points $p, q \in \hat{M}$ with $\dist(p, q) = \pi$.  We will show that then $\hat{M} = \calh(x)$.

	Since $\hat{M} \subset \calh(x)$ we know that $p, q \in \cals(x)$. Let $y\in \interior \calh(x)$ be arbitrary. Then there exists a unique $C^1$-geodesic $\gamma:[0, \pi]\rightarrow \bbs^{n+1}$ starting at $p$ and ending at $q$, such that $y \in \gamma((0, \pi)) \subset \interior \calh(x)$. From Lemma \ref{StereographicInball} applied to $p$ and $q$ we obtain that $\gamma((0, \pi)) \subset \hat{M}$. Hence $\interior\calh(x) \subset \hat{M}$ and we infer $\hat{M} = \ov{\hat{M}} = \calh(x)$.
\end{proof}

Hence we know that the weakly convex set $\hat{M}$ in Theorem \ref{RigidityResult} is a convex set in a hemisphere. Thus it remains to distinguish two cases: $\hat{M}$ does or does not contain a pair of antipodal points.

The proof of the following Lemma can be found in \cite[Chapter 3, Corollary 1]{FerreiraIN}. For the sake of completeness, we give a sketch of an elementary proof.
\Lem
\label{OpenHemisphereLemma}
	Let $n\geq 1$ and $\hat{M} \subset \bbs^{n+1}$ be a convex body in the sphere, which does not contain any antipodal points. Then $\hat{M}$ is contained in an open hemisphere.
\eLem
\begin{proof}
	Firstly, since $\hat{M}$ is closed and does not contain pairs of antipodal points, we have $r:= \max\{\dist(p, q): p, q\in \hat{M}\}<\pi$.

	We will prove the Lemma by induction on $n \in \bbn$. For $n=0$ the statement is obvious. Suppose we have proven the statement for $n-1 \in \bbn$.
	
	 Let $p_1, q_1 \in \hat{M}$ be two points with $\dist(p_1, q_1) = r$. Let $\epsilon := \frac{\pi - r}{2}$. We can assume without loss of generality, that 
	 \eq p_1 = (\sqrt{1-\epsilon^2}, 0, \ldots, 0, \epsilon)\eeq
	  and \eq q_1 = (-\sqrt{1-\epsilon^2}, 0, \ldots, 0, \epsilon).\eeq 	
	Furthermore $\hat{M}^1 := \hat{M} \cap \{x^1= 0\}$ is a closed subset of $\bbs^n \equiv \{0\}\times \bbs^n \subset \bbs^{n+1}$ and satisfies the requirements of the lemma for $m= n-1$. Hence by the inductive assumption, $\hat{M}^1$ is contained in an open hemisphere and we can assume after a rotation about the $x^1$-axis, that $\hat{M}^1 \subset \{x\in \bbs^{n+1}: x^{n+2}\geq \epsilon\}$.  
	
	Suppose there exists a point $z\in \hat{M}$ with $z^{n+2} = 0$.  In view of the observations made above there holds $0<|z^1|<1$, $\hat{z} = (z^2, \ldots, z^{n+1}) \neq 0$ and we can assume without loss of generality $-1<z^1<0$. Then 
	\begin{equation}
		\hat{M}\cap \{x^1\geq 0\} \subset \{x^{n+2}> 0\},
	\end{equation} 
	for otherwise there would exist a point $y\in \hat{M}\cap\{x^{n+2} = 0, x^1>0\}$ and either $\dist(y, z) = \pi$, which is excluded by the assumption of the lemma, or $\dist(y, z) <\pi$ and hence the geodesic segment $\hat\Gamma_{y, z}$ would be contained in $\hat{M} \cap\{y^{n+2} = 0\}$, which implies a contradiction to 
	\begin{equation}
		\hat{M}^1 \cap \{x^{n+2} = 0\} = \emptyset.
	\end{equation} 
	
	Now we rotate continuously in the positive $x^1$-direction such that $\hat{M}\cap\{x^1>0\}\subset \{x^{n+2}\geq 0\}$ and there exists $y\in \hat{M} \cap\{x^{n+2} = 0, x^1>0\}$.
	 Note that we still have $\hat{M}\cap\{x^1=0\} =\emptyset$. By the same reasoning as above, in the new coordinate system $\hat{M} \cap\{x^{n+2} =0, x^1 < 0\} = \emptyset$. This implies that $\hat{M}$ is contained in an open hemisphere. 
\end{proof}


\begin{proof}[Proof of Theorem \ref{RigidityResult}]
From Theorem \ref{WeaklyConvexTheorem} we obtain that $\hat{M}$ is a convex body in a hemisphere.  Hence if $\hat{M}$ contains a pair of antipodal points, then $\hat{M} = \calh(x_0)$. Otherwise it is contained in an open hemisphere in view of Lemma \ref{OpenHemisphereLemma}.

\end{proof}


For $C^{2}$-hypersurfaces we obtain a generalization of Theorem \ref{doCarmoWarner}, see Corollary \ref{C2Rigidity}. For the proof of this result, we need one further Lemma.

\Lem \label{hemisphere} 
Let $\hat{M}_{n}\subset \mc{H}({x_{n}})\subset \S^{n+1}$ be a sequence of sets, such that 
\eq \hat{M}_{n}\subset \inte\hat{M}_{n+1},\eeq
then there exists $x_{0}\in\S^{n+1},$ such that

\eq \label{h1} \hat{M}_{n}\subset \inte \mc{H}({x_{0}})\ \ \forall n\in\N.\eeq  
\eLem

\pf
A subsequence of points $x_{n_{k}}$ converges
to some $x_{0}\in\S^{n+1}.$ We claim that this $x_{0}$ is a point which satisfies (\ref{h1}). 
If this was not the case, then we use the monotonicity to derive the existence of $n_{0}\in\N$ with the property
\eq \hat{M}_{n_{0}}\cap \mc{H}(x_{0})^{c}\neq \emptyset.\eeq
Thus there exists a point
\eq y\in\hat{M}_{n},\ n\geq n_{0},\eeq
and $\e>0$ with the property
\eq \dist(y,\mc{H}(x_{0}))\geq \e.\eeq
This leads to a contradiction, since for large $k$ we have
\eq y\in \hat{M}_{n_{k}-1}\subset \inte \mc{H}(x_{n_{k}})\subset\inte \mc{H}(x_{n_{k}})\cup\inte \mc{H}(x_{0})\eeq
and the maximal distance of points in the latter set to $x_{0}$ converges to $\frac{\pi}{2}.$
\epf

\begin{proof}[Proof of Corollary \ref{C2Rigidity}]
Choose a differentiable normal vector field $\nu,$ such that the second fundamental form with respect to $-\nu$ is positive semi-definite and let $\mc{U}$ be a tubular neighborhood around $M$ with corresponding signed distance function 
$d\in C^{2}(\mc{U})$ and normal Gaussian coordinate system $(x^{\a}),$ 
compare \cite[Theorem 1.3.13]{cg:cp}. Note that $d=x^{0}$ and $\nabla d=\nu.$ According to \cite[Lemma 2.4.3]{cg:cp}, the second fundamental form with respect to $-\nu$ of the coordinate slices 
\eq \{x^{0}=\const\},\eeq
which can be seen as a solution to the flow
\eq \dot{x}=\nu,\ x(t,\xi)=(t,\xi),\eeq
evolves according to the evolution equation
\eq \dot{h}^{i}_{j}=-h^{k}_{j}h^{i}_{k}-\d^{i}_{j}.\eeq
Thus the principle curvatures of the slices are strictly decreasing, which implies, that the hypersurfaces
\eq M_{-t}=\{x^{0}=-t\}\eeq are strictly convex with positive definite second fundamental form. Consider the image $\mc{M}_{-t}$ under a suitable stereographic projection $\mc{P}$, which is a strictly convex $C^{2}$ hypersurface in $\R^{n+1}.$ For any $\d>0$ there exists $\e(\d)>0,$ such that for the convolution of the signed distance function, $d_{\e},$ there hold
\eq \mc{M}_{-t}^{\e}\equiv\left\{d_{\e}=-\(t+\frac{\d}{2}\)\right\}\subset \mc{P}(\{-(t+\d)<d<-t\}),\eeq
\eq \langle\nabla d_{\e},\nabla d\rangle\geq c> 0\eeq and
\eq \mc{M}_{-t}^{\e}\ \text{is strictly convex}.\eeq 
Using (\ref{Ster.Proj.A}) and the $C^2$-convergence of the convolution, those properties carry over to $M_{-t}^{\e}\equiv\mc{P}^{-1}(\mc{M}_{-t}^{\e})$, a hypersurface in $\S^{n+1},$ to which we may apply Theorem \ref{doCarmoWarner}.  Using the same construction, we obtain 
\eq M_{-\frac{t}{2}}^{\e_{1}},\ \text{where}\ \e_{1}=\e\(\frac{t}{2}\).\eeq
Thus we derive a sequence of smooth and strictly convex hypersurfaces 
\eq M_{n}\equiv M_{-\frac{t}{2^{n}}}^{\e_{n}},\ \text{where}\ \e_{n}=\e\(\frac{t}{2^{n}}\),\eeq
with the property
\eq {M}_{n}\subset\inte\hat{M}_{n+1}. \eeq Here we also used the generalized Jordan curve theorem, cf. \cite[Chapter IV, 19]{Bredon}.
Lemma \ref{hemisphere} implies that there exists $x_0\in \bbs^{n+1}$ such that $M\subset \calh(x_0)$ and Theorem \ref{WeaklyConvexTheorem} shows that $M$ bounds a convex body $\hat{M} = \overline{\bigcup_{n\in \bbn} \hat{M}_n}$ in a hemisphere, since $\hat{M}$ obviously satisfies the interior sphere condition with respect to $x_0$ at all points of $\hat{M}\cap\cals(x_0)$.	
\end{proof}

\section{The curvature flow and first estimates}

\vspace{0,1 cm}
\subsubsection*{Curvature functions}
Now we mention some elementary facts about curvature functions on a hypersurface. 
\label{Notation}


To derive the geometric inequalities, we will need some properties of the elementary symmetric polynomials. 

\Lem
\label{SymPolLemma}
Let $1\leq k \leq n$ be fixed.
	\begin{itemize}
		\item[(i)] We define the convex cone 
			\begin{equation}
				\Gamma_k = \{(\kappa_i) \in \bbr^n: H_1(\kappa_i) > 0, H_2(\kappa_i) > 0, \ldots, H_k(\kappa_i) > 0 \}.
			\end{equation}
			Then $H_k$ is strictly monotone on $\Gamma_k$ and $\Gamma_k$ is exactly the connected component of
			\begin{equation}
				\{(\kappa_i) \in \bbr^n: H_k(\kappa_i) > 0\}
			\end{equation}
			containing the positive cone.
		\item[(ii)] The $k$-th roots $\sigma_k = H_k^{\frac{1}{k}}$ are concave on $\Gamma_k$.
		\item[(iii)] For $ 1 < s < t < n$ and $\tilde{\sigma}_k = \Bigl(\frac{H_k}{{n \choose k}}\Bigr)^{\frac{1}{k}}$ there holds
			\begin{equation}
				\tilde{\sigma}_n \leq \tilde{\sigma}_t \leq \tilde{\sigma}_s \leq \tilde{\sigma}_1,
			\end{equation}
			where the principal curvatures have to lie in $\Gamma_n \equiv \Gamma_+$ for the first, in $\Gamma_t$ for the second and in $\Gamma_s$ for the third inequality.
		\item[(iv)] For fixed $i$, no summation over $i$, there holds
			\begin{equation}
				H_k = \frac{\partial H_{k+1}}{\partial \kappa_i} + \kappa_i \frac{\partial H_k}{\partial \kappa_i}.
			\end{equation}
	\end{itemize}
\eLem
\begin{proof}
	The convexity of the cone $\Gamma_k$ and (i) follows from \cite[Section 2]{HuiskSinestr}, (ii) follows from \cite[Thm.~3.2]{Gerhardt:/2015} ,(iii) from \cite[Lemma 15.12]{Lieberman} and (iv) follows directly from the definition of the $H_k$.
\end{proof}

A consequence of the preceding lemma is the following
\Lem
\label{DivFreeLemma}
	Let $N$ be a semi-Riemannian space of constant curvature, then the symmetric polynomials $F= H_k$, $1\leq k \leq n$, are divergence free for every admissible hypersurface $M$ of $N$. In case $k=2$ it suffices to assume that $N$ is an Einstein manifold.
\eLem
\begin{proof}
	The proof of the lemma can be found in \cite[Lemma 5.8]{GerhSurvey}. The proof consists of induction on $k$ and (iv) of Lemma \ref{SymPolLemma}.
\end{proof}

Now we state a well-known inequality for general curvature functions:
\Lem
\label{FHInequalityLemma}
	Let $F \in C^2(\Gamma_+)$ be a strictly monotone, concave (respectively convex) curvature function, positively homogeneous of degree 1 with $F(1,\ldots, 1) > 0$, then 
	\begin{equation}
	F(\kappa) \leq (\textnormal{respectively } \geq)\, \frac{F(1, \ldots, 1)}{n} H(\kappa).
	\end{equation}
	and
	\begin{equation}
		\sum_{i=1}^n F_i(\kappa) \geq (\textnormal{respectively } \leq)\, F(1, \ldots, 1),
	\end{equation}
	where $\kappa = (\kappa_k) \in \Gamma_+$.
\eLem
\begin{proof}
	See \cite[Lemma 2.2.19, Lemma 2.2.20]{cg:cp}.
\end{proof}

\subsubsection*{Evolution equations for the curvature flow}


The following evolution equations are valid for curvature flows in the sphere, for a derivation see \cite[Chapter 2]{cg:cp}.

\Lem (Evolution equations)
\label{EvolutionEquationsLemma}
\begin{equation} \frac{d}{dt}g_{ij}=-2\Phi h_{ij}.\end{equation}
\begin{align}	\frac{d}{dt}h_{ij}&-\Phi'F^{kl}h_{ij;kl}=\Phi'F^{kl}h_{rk}h^{r}_{l}h_{ij}-(\Phi'F+\Phi)h_{ri}h^{r}_{j}\\
						&+(\Phi'F+\Phi)g_{ij}-\Phi'F^{kl}g_{kl}h_{ij}+\Phi^{kl,rs}h_{kl;i}h_{rs;j}.\notag\end{align}		

\begin{align}
\label{EvCurvatures}
\dt h^i_j &- \Phi' F^{kl} h^i_{j;kl}  = \Phi' F^{kl}h^r_kh_{rl} h^i_j + (\Phi - \Phi'F) h^i_kh^k_j\\ 
&+ \Phi' F^{kl,rs}h_{kl;m}h_{rs;j} g^{mi} + \Phi'' F^iF_j \notag\\
& + ((\Phi+\Phi'F) \delta^i_j - \Phi' F^{kl}g_{kl} h^i_j). \notag\end{align}
\begin{align}
\label{EvPhi}
\dt \Phi &-\Phi' F^{kl} \Phi_{;kl} = \Phi' F^{kl}h_{kr}h^r_l\Phi + \Phi'  F^{kl}g_{kl} \Phi\end{align}

\eLem

\subsubsection*{Curvature estimates}

Since the flow reduces to a scalar parabolic equation, see for example \cite[Chapter 2.5]{cg:cp}, short-time existence is guaranteed and we know that the flow exists on a maximal time interval $[0, T^*)$ for some $0<T^*\leq \infty$. 

Firstly, we show that the flow exists only up to a finite time $T^*$.
\Lem
\label{FiniteTimeLemma}
	There holds
	\begin{equation}
		T^*<\infty.
	\end{equation}
\eLem
\begin{proof}
	Using Lemma \ref{FHInequalityLemma} we can compare the evolution of $-\Phi$, see \eqref{EvPhi}, with the solution of the ordinary differential equation 
	\begin{equation}
		\dot\varphi = pn \varphi^{\frac{2p+1}{p}}.
	\end{equation}
	This shows that $\underset{x\in M_t}{\sup} (-\Phi(x))$ becomes unbounded in finite time. 
\end{proof}

We show that the hypersurfaces remain strictly convex and the principal curvatures are uniformly bounded from above. For curvature functions with $F_{|\partial\Gamma_{+}} = 0$ the strict convexity of the flow hypersurfaces follows immediately. For concave and inverse concave curvature functions this follows from the following.

\Lem
\label{PinchingLemma}
	Suppose that $F$ is a curvature function as in \ref{MainAssumption} (i) with $p=1$.
Then there exists a constant $c>0,$ such that the principal curvatures of the flow $\k_{1},\dots,\k_{n}$ satisfy
\eq \k_{n}\leq c\k_{1}.\eeq
\eLem
\begin{proof}

Let $0<\e<\frac{1}{n},$ such that
\eq T_{ij}=h_{ij}-\e Hg_{ij}\eeq
is positive definite at $t=0.$ $T_{ij}$ satisfies
\begin{align}\begin{split}	T_{ij}&-\Phi'F^{kl}T_{ij;kl}=\Phi'F^{kl}h_{rk}h^{r}_{l}T_{ij}-\Phi'F^{kl}g_{kl}T_{ij}\\
						&+\e(\Phi'F-\Phi)\|A\|^{2}g_{ij}+2\e\Phi Hh_{ij}\\
						&+\Phi^{kl,rs}h_{kl;i}h_{rs;j}-\e\Phi^{kl,rs}h_{kl;m}{h_{rs}}^{m}g_{ij}\\
						&\equiv N_{ij}+\tilde{N}_{ij},\end{split}\end{align}
where we also used $\Phi'F+\Phi=0,$ due to $p=1.$ Here $\tilde{N}_{ij}$ denotes the terms involving derivatives of $h_{ij}.$ At a point $(t_{0},\xi_{0})$ let $\eta$ be a null-eigenvector of $T_{ij},$ i.e.
\eq h_{ij}\eta^{j}=\e H\eta_{i}\eeq
and $T_{ij}\geq 0$ elsewhere. There holds
\begin{align}\begin{split} N_{ij}\eta^{i}\eta^{j}&=2\e\Phi'F\|A\|^{2}\|\eta\|^{2}+2\e^{2}\Phi H^{2}\|\eta\|^{2}\\
									&\geq 2\e\Phi'F H^{2}\|\eta\|^{2}\(\frac 1n-\e\)>0.\end{split}\end{align}
To prove that $\tilde{N}_{ij}$ satisfies a modified null-eigenvector condition, note that $\Phi$ is, as a function of the principal curvatures, symmetric, monotone, concave and inverse concave. We apply \cite[Theorem 4.1]{AndrewsPinching} to obtain
\eq \tilde{N}_{ij}\eta^{i}\eta^{j}+2\sup_{\Gamma}\Phi^{kl}(2\Gamma^{r}_{k}(h_{ri;l}\eta^{i}-\e H_{l}\eta_{r}-\Gamma^{r}_{k}\Gamma^{q}_{l}(h_{rq}-\e H\delta_{rq})))\geq 0.\eeq
Those are exactly the requirements to apply Andrews' generalized maximum principle, \cite[Theorem 3.2]{AndrewsPinching} to conclude $T_{ij}>0$ for all $t\in[0,T^{*}).$

\end{proof}

Next, we derive upper bounds for the principal curvatures.

\Lem
\label{HBoundLemma}
	Let $c_1:= \sup_{M_0} H$. Then there holds 
	\begin{equation}
	\underset{t\in [0, T^*), x\in M_t}{\sup} H(x) \leq c_1.
	\end{equation}
\eLem
\begin{proof}
	
	Let $0<T<T^*$ be arbitrary. Suppose there exists $t_0\in (0, T]$ and $x_0\in M_{t_0}$ such that
	\begin{equation}
		\underset{t\in[0, t_0], x\in M_t}{\sup} H \leq H(x_0). 
	\end{equation}
	Then we obtain from the maximum principle and the concavity of $F$, that at $x_0$ there holds in view of \eqref{EvCurvatures}
	\begin{equation}
	\begin{split}
		0 &\leq  pF^{kl}h^r_kh_{rl} \frac{H}{F^{p+1}} - (p+1)\frac{|A|^2}{F^{p}} -  pF^{kl}g_{kl} \frac{H}{F^{p+1}} -(1-p) \frac{n}{F^p}\\
		& = \frac{p}{F^{p+1}} \( F^{kl}h^r_kh_{rl} H - F |A|^2\) -\frac{n+|A|^2}{F^p}\\
		&\hphantom{=} -\frac{p}{F^{p+1}}\(F^{kl}g_{kl} H - n F\).
	\end{split}
	\end{equation}
	However, we note that
	\begin{equation}
		\begin{split}
			F^{kl}h^r_kh_{rl} H - F\,|A|^2 &= \sum_{i,j} \left(f_i\kappa_i^2\kappa_j - f_i\kappa_i\kappa_j^2\right)\\
			&=\sum_{i<j}\kappa_i\kappa_j(\kappa_i-\kappa_j)\left(f_i-f_j\right)\leq 0,\\
		\end{split}
	\end{equation}
		since for concave curvature functions there holds $f_i \geq f_j$ for $i<j$. Furthermore, in view of Lemma \ref{FHInequalityLemma} we obtain
		\begin{equation}
			F^{kl}g_{kl} H - nF \geq 0.
		\end{equation}
	Hence we obtain a contradiction.
\end{proof}

\section{Elementary flow properties and further curvature estimates}

\Lem \label{monotonicity}
Let $M_{t},$ $0\leq t<T^{*},$ be a flow hypersurface and $\hat{M}_{t}$ be the enclosed convex body, cf. Corollary \ref{C2Rigidity}.  Then those convex bodies are strictly monotonically ordered, i.e.
\eq s<t\Ra \hat{M}_{s}\subset\mrm{int}\ \hat{M}_{t}.\eeq
\eLem

\pf
Let $0\leq s<T^{*}.$ Then the flow hypersurface $M_{s}$ is strictly convex, as was shown in the previous section. Thus, using Corollary \ref{C2Rigidity}, we first conclude that $M_{s}$ does indeed enclose a convex body $\hat{M}_{s}$ and that this body has to lie compactly in an open hemisphere $\inte\mc{H}(\~{x}_{s}).$ Choose $x_{s}\in \mathrm{int}\ \hat{M}_{s}.$ 
Thus $M_{s},$ which in particular is starshaped with respect to $x_{s},$ can be written as a graph over $\mc{S}({x_{s}}),$
\eq M_{s}=\graph u(s,\mc{S}({x_{s}})), \eeq
and thus there is an $\e>0,$ such that for all $s\leq t<s+\e$ the hypersurfaces $M_{t}$ may be written as a graph over $S({x_{s}}),$ compare \cite[Theorem 2.5.19]{cg:cp}. In these coordinates $u$ locally satisfies the scalar flow equation
\eq \frac{\del u}{\del t}=\frac{v}{F^{p}},\eeq
cf. \cite[p. 98-99]{cg:cp}, 
and thus the function $u$ is strictly increasing for fixed $x\in\mc{S}({x_{s}}).$
Since 
\eq \hat{M}_{s}=\{(r,(x^{i}))\in \R\times \mc{S}({x_{s}})\cn 0\leq r\leq u(s,(x^{i}))\},\eeq 
where $(r,(x^{i}))$ describe the corresponding geodesic polar coordinates around $x_{s},$
the claim follows.
\epf

\Prop \label{limitsurface}
\label{C1Conv}
There is a uniquely determined limit surface $M_{T^{*}},$ which can be written as a graph in geodesic polar coordinates,
\eq M_{T^{*}}=\graph u(T^{*},\mc{S}({y_{0}})),\eeq 
where $y_{0}\in\inte \hat{M}_{T^{*}}.$
Furthermore there holds
\eq u(t,\cdot)\ra u(T^{*},\cdot)\ \ \mrm{in}\ C^{1,\b}(\mc{S}({y_{0}}))\ \ \forall 0\leq \b<1.\eeq
\eProp

\pf
In those geodesic polar coordinates the metric of $\S^{n+1}$ is given by
\eq d\-{s}^{2}=dr^{2}+\sin^{2}r\s_{ij}dx^{i}dx^{j}, \ 0<r<\pi. \eeq 
Let $y_{0}\in\inte\hat{M}_{0}$ and $\hat{y}_{0}$ denote the antipodal point of $y_{0},$ then by Lemma \ref{hemisphere} we know that
\eq \hat{M}_{t}\subset K\Subset \S^{n+1}\backslash\{\hat{y}_{0}\}\ \ \forall 0\leq t<T^{*}.\eeq
Thus we have a uniform parametrization of the flow hypersurfaces as graphs over $\mc{S}({y_{0}}),$
\eq M_{t}=\graph u(t,\cdot).\eeq
The quantity $v^{2}=1+\sin^{-2}u\s^{ij}u_{i}u_{j}$ is bounded by convexity, see \cite[Theorem 2.7.10]{cg:cp}.
The second fundamental form of a graph hypersurface satisfies 
\eq h^{i}_{j}=\frac{\dot{\vt}}{v\vt}\d^{i}_{j}+\frac{\dot{\vt}}{v^{3}\vt^{3}}u^{i}u_{j}-\frac{\~{\s}^{ik}}{v\vt^{2}}u_{kj},\eeq
where $\vt=\sin u$ and $\~{\s}^{ik}$ is the inverse of 
\eq \~{\s}_{ik}=\p_{i}\p_{k}+\s_{ik},\ \ \p=\int_{r_{0}}^{u}\sin^{-1}(s) ds,\eeq
cf. \cite[(3.82)]{js:ImpfHyp}. Here covariant differentiation and index raising is performed with respect to the metric $\s_{ij}$
and by Lemma \ref{HBoundLemma} we obtain uniform $C^{2}$ estimates for the $u(t,\cdot).$ This gives the existence of a convergent subsequence with uniquely determined $C^{1}$ limit $u(T^{*},\cdot)$, using monotonicity.

\epf
\Prop \label{Fpositive}
Suppose that on some time interval $[s,t]\subset (0,T^{*})$ there is a point $x_{0}\in \mrm{int}\ \hat{M}_{s},$ such that for a common parametrization of the surfaces $M_{\tau},$
\eq M_{\tau}=\graph u(\tau,\mc{S}({x_{0}})),\eeq
there is a constant $\e>0$ satisfying  $u\leq\frac{\pi}{2}-\e,$
then the curvature function of those hypersurfaces is uniformly positive,
\eq  F\geq \~{c}>0,\eeq where $\~{c}$ depends on $\e.$

\eProp

\pf
Using \cite[Lemma 3.3.2]{cg:cp} and (\ref{EvPhi}) we deduce the evolution equations for $-\Phi$ and $u$ to be
\eq \frac{d}{dt}(-\Phi)-\Phi'F^{ij}(-\Phi)_{ij}=\Phi'F^{ij}h_{ik}h^{k}_{j}(-\Phi)+\Phi'F^{ij}g_{ij}(-\Phi) \eeq
and 
\eq \frac{d}{dt}u-\Phi'F^{ij}u_{ij}=(p^{-1}+1)\Phi'Fv^{-1}-\Phi'F^{ij}\-{h}_{ij}.\eeq 

Set \eq w=\log\left(-\Phi\right)+f(u),\eeq
where $f$ will be specified later.
Then
\begin{align} \begin{split} \frac{d}{dt}w-\Phi'F^{ij}w_{ij}&=\Phi'F^{ij}h_{ik}h^{k}_{j}+\Phi'F^{ij}g_{ij}+\Phi'F^{ij}\frac{\Phi_{i}}{\Phi}\frac{\Phi_{j}}{\Phi}\\
				&\hp{=}+(p^{-1}+1)f'\Phi'Fv^{-1}-f'\Phi'F^{ij}\-{h}_{ij}-f''\Phi'F^{ij}u_{i}u_{j}.\end{split}\end{align}
We want to bound the function $w.$ Thus suppose without loss of generality, that 
\eq \sup\limits_{(\tau,\xi)\in(s,t]\times M}w(\tau,\xi)=w(t_{0},\xi_{0}).\eeq
Then at this point we have
\eq \frac{\Phi_{i}}{\Phi}=-f'(u)u_{i}\eeq
and thus, also using $\-{h}_{ij}=\frac{\-{H}}{n}\-{g}_{ij}=\frac{\-{H}}{n}g_{ij}-\frac{\-{H}}{n}u_{i}u_{j},$ we find at $(t_{0},\xi_{0})$
\begin{align}
\begin{split} 
0&\leq \Phi'F^{ij}h_{ik}h^{k}_{j}+\(1-f'\frac{\-{H}}{n}\)\Phi'F^{ij}g_{ij}+(p^{-1}+1)f'\Phi'Fv^{-1}\\
					&\hp{=}+\((f')^{2}+f'\frac{\-{H}}{n}-f''\)\Phi'F^{ij}u_{i}u_{j}.\end{split}\end{align}
Now define the function 
\eq f(u)=-\log\(\cos u-\cos\(\frac{\pi}{2}-\frac{\e}{2}\)\)=-\log(\cos u-c),\eeq
where $c:=\cos\(\frac\pi2 - \frac\epsilon 2\)$, such that  by assumption $f$ is strictly positive and uniformly bounded for $\tau\in[s,t]$.
We have
\eq f'=\frac{\sin u}{\cos u-c}\eeq and
\eq f''=\frac{1}{(\cos u-c)^{2}}-\frac{c\cdot \cos u}{(\cos u-c)^{2}}\eeq
and thus in view of \eqref{HLevelHyp}
\eq f'\frac{\-{H}}{n}=\frac{\cos u}{\cos u-c}>1+\d,\ \d=\d(\e)\eeq and
\eq (f')^{2}+f'\frac{\-{H}}{n}-f''=0.\eeq
Since $F^{ij}h_{ik}h^k_j \leq F H$ due to the convexity of the hypersurfaces, we conclude at $(t_{0},\xi_{0})$
\eq 
\begin{split}
	0&\leq \Phi'F^{ij}h_{ik}h^{k}_{j} +(p^{-1}+1)f'\Phi'Fv^{-1} -n\d\Phi'\\
	&\leq \frac{p}{F^p} H + (p+1) \frac{f'}{vF^{p}} - n\delta \frac{p}{F^{p+1}}. 
\end{split}
\eeq
Supposing that $w(t_{0},\xi_{0})$ is very large, $-\Phi$ must also be very large, which leads to a contradiction, since $H$ is bounded by Lemma \ref{HBoundLemma}. Hence $w,$ and thus also $-\Phi,$ must be bounded. 
\epf

Now we characterize $T^{*},$ where in particular we show, that the flow exists as long as the hypersurfaces remain strictly convex.

\Prop \label{Fzero}
Suppose that on some interval $[s,t)\subset [0,T^{*})$ we have
\eq F(\tau,\xi)\geq \e>0\ \ \forall(\tau,\xi)\in[s,t)\times\S^{n},\eeq
then there holds
\eq T^{*}>t.\eeq
\eProp

\pf
Without loss of generality suppose that all $M_{\tau}, s\leq \tau< t,$ are uniformly parametrized over $\mc{S}({y_{0}}).$
The second fundamental form has, with respect to the corresponding spherical metric, the form
\eq h^{i}_{j}=\frac{\dot{\vt}}{v\vt}\d^{i}_{j}+\frac{\dot{\vt}}{v^{3}\vt^{3}}u^{i}u_{j}-\frac{\~{\s}^{ik}}{v\vt^{2}}u_{kj}.\eeq
The scalar flow equation is
\eq \frac{\del}{\del t}u=\frac{v}{F^{p}(h^{i}_{j})}\equiv G(x, u,Du,D^{2}u).\eeq
Thus
\eq \frac{\del G}{\del u_{ij}}=-\frac{pv}{F^{p+1}}F^{k}_{l}\frac{\del h^{l}_{k}}{\del u_{ij}}=\frac{pv}{F^{p+1}}F^{k}_{l}\frac{\~{\s}^{lm}}{v\vt^{2}}\d^{i}_{m}\d^{j}_{k}=\frac{p}{\vt^{2}F^{p+1}}F^{j}_{l}\~{\s}^{li},\eeq
where the latter is uniformly positive definite by assumption, as well as by Lemma \ref{PinchingLemma} or $F|_{\partial\Gamma_+} = 0$. Then, using Krylov-Safonov and Schauder, we obtain uniform $C^{4,\alpha}$ 
estimates on $[s,t)$ and thus the flow extends beyond $t.$
\epf

\section{Convergence to an equator}
In view of Lemma \ref{monotonicity} we know that there exists $x_0 \in \bbs^{n+1}$, such that $\hat{M_t} \subset \interior \calh(x_0)$. Now we want to show, that the limit hypersurface $M_{T^{*}}$ is equal to $\mc{S}(x_{0}).$ We first view the hypersurfaces as embedded in Euclidean space using the conformally flat parametrization of $\S^{n+1}.$

\Lem
\label{InteriorSphereLemma}
The enclosed, weakly convex body of $M_{T^{*}},$ $\hat{M}_{T^{*}},$ satisfies a uniform interior sphere condition.
\eLem

\pf
We will denote $\calp_{x_0}\left(M_t\right)$ by $\mc{M}_t$ for $t\in [0, T^*]$.
Since all $M_{t}$ range within distance less than $\frac{\pi}{2}$ around $x_{0},$ the metrics $\-{g}_{\a\b}$ and $\hat{g}_{\a\b}$ are uniformly equivalent on the set of consideration. Thus we also obtain the $C^{1}$-convergence of $\mc{M}_{t}\ra \mc{M}_{T^{*}}.$
Let $x\in\mc{M}_{T^{*}}$ be arbitrary and $t_{n}, x_{n}\in\mc{M}_{t_{n}}$ be sequences such that
\eq t_{n}\ra T^{*},\ x_{n}\ra x.\eeq
By Remark \ref{StereographicProjection} we obtain a sequence of inballs with center $y_{n}\in \inte\hat{\mc{M}}_{t_{n}}$ and uniform radius $R,$ such that
\eq B_{R}(y_{n})\subset \inte\hat{\mc{M}}_{t_{n}}\ \wedge\ \del B_{R}(y_{n})\cap \mc{M}_{t_{n}}=\{x_{n}\}.\eeq
Without loss of generality we have
\eq y_{n}\ra y\in\inte\hat{\mc{M}}_{T^{*}}.\eeq
First of all, let $z\in B_{R}(y).$ By the triangle inequality for large $n$ there holds \eq z\in B_{R}(y_{n})\subset \inte\hat{\mc{M}}_{t_{n}}\subset\inte\hat{\mc{M}}_{T^{*}}.\eeq
Thus \eq B_{R}(y)\subset \inte\hat{\mc{M}}_{T^{*}}.\eeq
There holds
\eq \dist(x,y)\leq \dist(x,x_{n})+\dist(x_{n},y_{n})+\dist(y_{n},y)\ra R\eeq
and thus \eq \dist(x,y)=R,\eeq 
since $x\in\mc{M}_{T^{*}}.$ We summarize:
\eq B_{R}(y)\subset\inte\hat{\mc{M}}_{T^{*}}\ \wedge\ \{x\}\subset \del B_{R}(y)\cap\mc{M}_{T^{*}}.\eeq
If we now choose an inball \eq B_{\frac{R}{2}}(\~{y})\subset B_{R}(y)\eeq with the property
\eq \del B_{\frac{R}{2}}(\~{y})\cap \del B_{R}(y)=\{x\},\eeq
we obtain the desired uniform interior sphere condition with radius $\frac{R}{2}.$ 
\epf

\Lem
\label{GeodesicConvexLemma}
	$\hat{M}_{T^*}\subset \calh(x_0)$ is a convex body in a hemisphere.
\eLem
\begin{proof}
	Let $y_0 \in \inte\hat{M}_0$ be arbitrary but fixed. In view of Lemma \ref{monotonicity} we know that the weakly convex body of $M_{T^*}$ with respect to  $y_0$ can be described as
\begin{equation}
\label{ConvexBody}
	\hat{M}_{T^*} = \overline{\bigcup_{t\in [0, T^*)} \hat{M}_t}.
\end{equation}
	In view of the monotonicity, see Lemma \ref{monotonicity} and Lemma \ref{hemisphere}, we obtain that $\bigcup_{t\in [0, T^*)} \hat{M}_t$ is a convex set in a hemisphere. Lemma \ref{InteriorSphereLemma} and Theorem \ref{WeaklyConvexTheorem} imply that $\hat{M}_{T^*}\subset \bbs^{n+1}$ is a convex body in a hemisphere.
\end{proof}
	
The lemmata of this section show, that $\hat{M}_{T^*}$ satisfies the requirements of Theorem \ref{RigidityResult}. Finally we will show that $\hat{M}_{T^*}$ is not contained in an open hemisphere:
\Lem
\label{LimitClosedHemisphere}
There is no hemisphere $\mc{H}(y_{0}),$ such that
\eq \hat{M}_{T^{*}}\equiv\overline{\bigcup_{t\in[0,T^{*})}\hat{M_{t}}}\subset \inte\mc{H}(y_{0}).\eeq 
\eLem

\pf
Suppose contrarily that there existed such a hemisphere. We may assume without loss of generality, that $y_{0}\in\inte\hat{M}_{T^{*}},$ compare \cite{Santalo:07/1946}.
Thus we are now in the situation, that we may parametrize the surfaces 
\eq M_{t},\ T<t<T^{*},\eeq
uniformly as a graphs over the same equator and may apply Propositions \ref{Fpositive} and \ref{Fzero} to conclude that the flow would exist longer than $T^{*}.$
\epf

Hence we obtain:
\Thm
\label{Equator}
	There exists $x_0\in \bbs^{n+1}$ such that $M_{T^*} = \cals(x_0)$.
\eThm

\section{Geometric inequalities}

In this section we want to deduce Alexandrov-Fenchel type inequalities from the convergence of the flow to an equator. We use the inverse mean curvature flow. The mixed volumes of a hypersurface $M$ in the sphere $\bbs^{n+1}$ are defined for $k\in \{0, \ldots, n\}$ by
\begin{equation}
	V_k(M) := \int_{M} \tilde{H}_k d\mu,
\end{equation}
where $\tilde{H}_k := \frac{H_k}{{n\choose k}}$ are the normalized elementary symmetric polynomials.

Firstly, we note that a geodesic ball $B_\rho$ of radius $0<\rho\leq \frac{\pi}{2}$ satisfies for $k\in \{0, \ldots, n\}$
\begin{equation}
	V_k(\del B_\rho) = \omega_n \cos^k\rho \sin^{n-k}\rho,
\end{equation} 
where $\omega_n$ is the volume of $\bbs^n$.
For the sake of brevity, we define $\tilde{V}_k := \frac{V_k}{\omega_n}$. For geodesic spheres, it is easy to obtain a relation between different mixed volumes. For example there holds
\begin{equation}
\label{GeomIneqSpheresI}
	(\tilde{V}_1)^2 = (\tilde{V}_0)^{2\left(\frac{n-1}{n}\right)} - (\tilde{V}_0)^2,
\end{equation}	
and
\begin{equation}
\label{GeomIneqSpheresII}
	\tilde{V}_2 = \left(\tilde{V}_0\right)^{\frac{n-2}{n}} - \tilde{V}_0.
\end{equation}
Furthermore, we will establish geometric inequalities between certain quermassintegrals. 

In $\bbs^{n+1}$ we have the following definition of the quermassintegrals, compare \cite{SolanesInt} (also for a more detailed definition of the measure $dL_k$):
\begin{defn}
	Let $\hat{M} \subset \bbs^{n+1}$ be a compact domain. For $k\in \{1, \ldots, n\}$ set
	\begin{equation}
		W_k(\hat{M}) = \frac{(n+1-k)\omega_{k-1}\cdots \omega_0}{(n+1)\omega_{n-1}\cdots\omega_{n-k}} \int_{\mathcal{L}_k}\chi(L_k\cap\hat{M})dL_k,
	\end{equation}
	where $\mathcal{L}_k$ is the space of $k$-dimensional totally geodesic subspaces $L$ in $\bbs^{n+1}$, $dL_k$ is the natural measure on $\mathcal{L}_k$ and $\chi$ is the Euler characteristic. Furthermore set
	\begin{equation}
		W_0(\hat{M}) = \Volu (\hat{M})
	\end{equation}
	and
	\begin{equation}
		W_{n+1}(\hat{M}) = \frac{\omega_n}{n+1} \chi(\hat{M}).
	\end{equation}
\end{defn}

In Euclidean space the quermassintegrals differ only by constants with respect to the corresponding curvature integrals. This relation is more complicated in curved spaces, however, we still have the following relation between the curvature integrals and the quermassintegrals in the space $\bbs^{n+1}$, see for example \cite[Proposition 7, Corollary 8]{SolanesInt} for a proof of this relation.

\begin{lemma}
\label{QuerVolLemma}
	If $\hat{M}\subset \bbs^{n+1}$ is a compact domain with $C^2$-boundary. Then there holds $W_1(\hat{M}) = \frac{1}{n+1} V_0(\partial M)$ and for $k\in \{1, \ldots, n\}$ there holds
	\begin{equation}
	\label{QuerVolEqI}
		V_k(\partial \hat{M}) = (n+1)\(W_{k+1}(\hat{M}) - \frac{k}{n+2-k} W_{k-1}(\hat{M})\).
	\end{equation} 
	Furthermore for $k\in \bbn$ with $2k+1\leq n$ we have
	\begin{equation}
	\label{QuerVolEqII}
		W_{2k+1}(\hat{M}) = \frac{1}{n+1} \sum_{i=0}^k\frac{(2k)!!(n-2k)!!}{(2k-2i)!!(n-2k+2i)!!} V_{2k-2i}(\partial\hat{M}).
	\end{equation}
\end{lemma}

We also need the following evolution equations. These evolution equations can be computed by an induction argument using Lemma \ref{EvolutionEquationsLemma} and Lemma \ref{QuerVolLemma}, see also the proof of \cite[Proposition 3.1]{WangXia}.

\begin{align}
\label{MixedVolEq}
&\frac{d}{dt} \int_{M_{t}} H_k d\mu_t = \int_{M_t} \left((k+1) \frac{H_{k+1}}{H} - (n+1-k) \frac{H_{k-1}}{H} \right) d\mu_t. \\
\label{QuermassintegralEq}
&\frac{d}{dt} W_k(M_t) = \frac{n+1-k}{n+1}\int_{M_t} \frac{\tilde{H}_k}{H} d\mu_t.
\end{align}

Finally, let us state a decay lemma for the inverse curvature flow.

\Lem\label{CurvatureDecayLemma}
For all $1\leq q<\infty$ there holds
\eq \lim\limits_{t\ra T^{*}}\int_{M_{t}}H^{q}=0.\eeq
\eLem

\pf
By the previous results we know
\eq u\ra c\eeq 
in the $C^{1}$ norm, where we use a graph representation as in Proposition \ref{limitsurface}. The second fundamental form is
\eq h_{ij}v^{-1}=-u_{ij}+\-{h}_{ij},\eeq
where covariant differentiation is performed with respect to the induced metric.
Thus
\eq \int_{M_{t}}Hv^{-1}=\int_{M_{t}}\-{h}_{ij}g^{ij}\eeq
and
\eq \int_{M_{t}}H=\int_{M_{t}}H(1-v^{-1})+\int_{M_{t}}\(\-{H}-\frac{\-{H}}{n}\|Du\|^{2}\).\eeq
Since
\eq g_{ij}=u_{i}u_{j}+\sin^{2}u\s_{ij},\eeq the volume element is uniformly bounded and thus the right hand side converges to $0.$ The other $L^{q}$ norms converge to $0$ by interpolation.

\epf

This leads to our first geometric inequality:

\Thm
\label{GeomIneqThmI}
	Let $M\subset \bbs^{n+1}$ be an embedded, closed, connected and convex $C^2$-hypersurface. Then we have the inequality
	\begin{equation}
		(\tilde{V}_1(M))^2 \geq (\tilde{V}_0(M))^{2\left(\frac{n-1}{n}\right)} - (\tilde{V}_0(M))^2,
	\end{equation}
	and equality holds if and only if $M$ is a geodesic sphere.
\eThm
\begin{proof}
	First of all, we can assume that the hypersurface is smooth and strictly convex, since otherwise we can use convolutions as in the proof of Corollary \ref{C2Rigidity} to obtain a sequence of approximating smooth, strictly convex hypersurfaces converging in $C^2$ to $M$.
	
	We consider the flow of the initial hypersurface $M$ by the inverse mean curvature. Let $M_t$, $t\in [0, T^*],$ be the level hypersurfaces of the flow, where we know that $M_{T^*}$ is a geodesic sphere with radius $\frac{\pi}{2}$ and the graphs over a geodesic sphere $M_t = $ graph$|_{\bbs^n} u(t, \cdot)$ converge in $C^1$ to $u(T^*, \cdot) \equiv \frac{\pi}{2}$. Furthermore we know that the mean curvature of the level hypersurfaces converges almost everywhere to $0$ by Lemma \ref{CurvatureDecayLemma}.
	 
	For $t\in [0, T^*]$ we define 
	\begin{equation}
		\phi(t) := \frac{(V_1(M_t))^2}{(V_0(M_t))^{2\left(\frac{n-1}{n}\right)}} + (V_0(M_t))^{\frac{2}{n}}.
	\end{equation}
 Then $\phi(t) \rightarrow (V_0(B_{\frac{\pi}{2}}))^{\frac2n}$ for $t\rightarrow T^*$, since $V_0(M_t)\rightarrow V_0(M_{T^*})$ for $t\rightarrow T^*$ in view of the $C^1$-convergence and $V_1(M_t) \rightarrow 0$ for $t\rightarrow T^*$ in view of the convergence of the mean curvature to zero almost everywhere.
	
	Hence if we can show that $\phi$ is monotonically non-increasing, we obtain
	\begin{equation}
		\frac{(V_1(M))^2}{(V_0(M))^{2\left(\frac{n-1}{n}\right)}} + (V_0(M))^{\frac{2}{n}} = \phi(0) \geq \phi(T^*) = \omega_n^{\frac2n}.
	\end{equation}
	This implies the geometric inequality.
	
	We have in view of \eqref{MixedVolEq} and (iii) from Lemma \ref{SymPolLemma}:
	\begin{equation}
	\begin{split}
		(V_0&(M_t))^{2\left(\frac{n-1}{n}\right) +1} \frac{d\phi}{dt} = 2V_0(M_t) \int_{M_t} \frac{H}{n}d\mu_t\int_{M_t}\left(\frac{2H_2}{nH} - \frac{1}{H}\right) d\mu_t \\
		&-2\frac{n-1}{n^3}\left(\int_{M_t} H d\mu_t\right)^2 \int_{M_t} 1~ d\mu_t + \frac2n (V_0(M_t))^2 \int_{M_t} 1~ d\mu_t\\
		&\leq \frac{2}{n}V_0(M_t)\left((V_0(M_t))^2 - \int_{M_t} Hd\mu_t\, \int_{M_t}\frac1H d\mu_t\right) \leq 0.   
	\end{split}
	\end{equation}
	In the last inequality we have used the Cauchy-Schwarz inequality and hence the inequality is strict unless $H$ is constant on $M_t$. Hence we obtain that $\phi$ is monotonically decreasing unless $M_t$ is a geodesic sphere.
\end{proof}

We can also prove a geometric inequality relating $\tilde{V}_2$ and $\tilde{V}_0$.

\Thm
\label{GeomIneqThmII}
	Let $n\geq 3$ and $M\subset \bbs^{n+1}$ be an embedded, closed, connected and convex $C^2$-hypersurface. Then we have the inequality
	\begin{equation}
	\label{GeomIneqEqII}
		\tilde{V}_2(M) \geq \(\tilde{V}_0(M)\)^{\frac{n-2}{n}} - \tilde{V}_0(M),
	\end{equation}
	and equality holds if and only if $M$ is a geodesic sphere.
\eThm
\begin{proof}
	Again we assume $M$ to be smooth and strictly convex and we use the same curvature flow as in the proof of Theorem \ref{GeomIneqThmI}. For $t\in [0, T^*]$ we define
	\begin{equation}
		\phi(t) := \frac{\(V_2(M_t)\) + \(V_0(M_t)\)}{\(V_0(M_t)\)^{\frac{n-2}{n}}}.
	\end{equation}
	Again we have to show that $\phi$ is monotonically non-increasing to obtain the geometric inequality \eqref{GeomIneqEqII}.
	We have in view of \eqref{MixedVolEq}:
	\begin{equation}
	\begin{split}
		\(V_0(M_t)\)&^{\frac{n-2}{n} + 1} \frac{d\phi}{dt} = V_0(M_t) \int_{M_t} \frac{3 H_3}{Ê{n\choose 2} H} - V_0(M_t) \int_{M_t} \(1-\frac{(n-1)}{{n\choose 2} }\)  d\mu_t\\
		& - \frac{n-2}{n} \int_{M_t} 1~ d\mu_t \(\int_{M_t} \tilde{H}_2 d\mu_t + V_0(M_t)\)\\
		&= \frac{V_0(M_t)}{ {n\choose 2}} \int_{M_t} \(\frac{3H_3}{H} - \frac{n-2}{n} H_2 \) d\mu_t \leq 0,
	\end{split}
	\end{equation}
	where the last inequality follows from (iii) from Lemma \ref{SymPolLemma}. Again, this inequality is strict unless the hypersurface is totally umbilic, which implies that $\phi$ is only stationary for geodesic spheres.
\end{proof}

\Thm
\label{GeomIneqThmIII}
	Let $M\subset \bbs^{n+1}$ be an embedded, closed, connected and convex $C^2$-hypersurface. Let $k\in \bbn_+$ with $2k+1\leq n$ and let $\hat{M}$ be the convex body enclosed by $M$. Then we have the inequality
	\begin{equation}
	\label{GeomIneqEqIII2}
		W_{2k+1}(\hat{M}) \geq \frac{\omega_n}{n+1}Ê\sum_{i=0}^k(-1)^i\frac{n-2k}{n-2k+2i} {k\choose i} \left(\frac{n+1}{\omega_n}W_1(\hat{M})\right)^{\frac{n-2k+2i}{n}}.
	\end{equation}	
	and equality holds if and only if $M$ is a geodesic sphere.
\eThm
\begin{proof}
	For the sake of brevity, let us denote the right hand side of the equation \eqref{GeomIneqEqIII2} by  $A_k(\hat{M})$.
	
	We prove the result by induction on $k$. For $k=1$ the result follows from Theorem \ref{GeomIneqThmII} and Lemma \ref{QuerVolLemma}. Let us assume we have proved the result for $k-1\in \bbn$ with $2k+1\leq n$. 
	
	We assume the hypersurface to be smooth and strictly convex and use the same curvature flow as in the proof of Theorem \ref{GeomIneqThmI}. For $t\in [0, T^*]$ we define
	\begin{equation}
		\phi(t) := \frac{W_{2k+1}(\hat{M}_t) - A_k(\hat{M}_t)}{W_1(\hat{M}_t)^{\frac{n-2k}{n}}}.
	\end{equation}
	From now on we will drop the arguments and write simply $W_k$ instead of $W_k(\hat{M})$.
	We will show that $\phi$ is monotonically non-increasing. We obtain from \eqref{QuermassintegralEq} and Lemma \ref{QuerVolLemma}:
	\begin{equation}
	\begin{split}
		\frac{d}{dt} \phi &\leq \frac{(n-2k)}{n(n+1)} \frac{\int_{M_t}\tilde{H}_{2k} d\mu_t - (n+1)W_{2k+1}}{W_1^{\frac{n-2k}{n}}} -\frac{d}{dt} \(\frac{A_k}{W_1^{\frac{n-2k}{n}}}\) \\
		& = -\frac{(n-2k)2k}{n(n+2-2k)}Ê\frac{W_{2k-1}}{W_1^{\frac{n-2k}{n}} }-\frac{d}{dt} \(\frac{A_k}{W_1^{\frac{n-2k}{n}}}\).
	\end{split}
	\end{equation}
	Note that the inequality at time $t$ is an equality if and only $M_t$ is totally umbilic. 
	Using the induction hypothesis we obtain after a simple calculation 
	\begin{equation}
		\frac{d}{dt} \phi \leq 0,
	\end{equation}
	again with equality if and only if $M_t$ is totally umbilic. Employing \eqref{QuerVolEqII} and the decay Lemma, we obtain 
	\begin{equation}
	\begin{split}
		W_{2k+1}(\hat{M}_t) - A_k(\hat{M}_t) &\rightarrow \frac{\omega_n}{n+1}\frac{(2k)!!(n-2k)!!}{n!!}\\
		& - \frac{\omega_n}{n+1}\sum_{i=0}^k (-1)^i \frac{n-2k}{n-2k+2i} {k\choose i}
	\end{split}
	\end{equation}
	for $t\rightarrow T^*$. However, the right hand side equals zero, as a proof by induction shows. Hence we obtain the inequality \eqref{GeomIneqEqIII2} and equality holds if and only if $M$ is a geodesic sphere.
\end{proof}

\textbf{Acknowledgements.}
We would like to thank Prof.~Oliver Schn\"urer for an interesting discussion concerning the interior sphere condition. We would like to thank Prof.~Claus Gerhardt for pointing out, that the class of curvature functions for which Theorem \ref{mainthm1} is true had to be modified. Furthermore we would like to thank Prof.~Rolf Schneider for valuable comments on isoperimetric problems for the mixed volumes.

\end{document}